\documentclass[12pt,utf8,a4paper,french,upmath]{aclart}
\usepackage{aclraccourcis,hyperref}

\def\Inf{\qopname\relax m{inf\vphantom p}}
\def\Sup{\qopname\relax m{sup\vphantom f}}

\let\mathmat\mathit
\title{Potentiel et rationalité}
\author{Antoine Chambert-Loir}
\address{%
Université de Paris, Sorbonne Université, CNRS, Institut de Mathématiques de Jussieu-Paris Rive Gauche, IMJ-PRG, F-75013, Paris, France}
\email{antoine.chambert-loir@u-paris.fr}

\author{Camille Noûs}
\address{Laboratoire Cogitamus}

\begin{abstract}
Nous étendons aux courbes de genre arbitraire le théorème de rationalité 
de Cantor, lui-même une extension de théorèmes de Borel, Pólya, Dwork, 
Bertrandias et Robinson.
La démonstration s'effectue en deux étapes. La première est un critère
d'algébricité, démontré par une méthode d'approximation diophantienne.
La seconde repose sur le théorème de l'indice de Hodge en théorie
d'Arakelov.
\end{abstract}
\begin{altabstract}
We extend to algebraic curves of arbitrary genus the rationality
theorem of Cantor, itself an extension of theorems of Borel, Pólya, Dwork,
Bertrandias and Robinson.
The proof runs in two steps. The first step is an algebraicity criterion,
which is proved using a method of diophantine approximation.
The second step relies on the Hodge index theorem in Arakelov geometry.
\end{altabstract}
\begin{document}
\maketitle

{\raggedleft
\itshape À Yuri, pour son soixantième anniversaire\par}

\section{Introduction}

Le théorème de \citet{Borel-1894} affirme 
qu'une série formelle à coefficients entiers 
est le développement de Taylor d'une fraction rationnelle
pourvu qu'elle définisse une fonction méromorphe dans un disque
de rayon~$>1$. Il a ensuite été étendu dans plusieurs directions. 
\citet{Polya-1928} a remplacé le disque par un voisinage
d'un compact de diamètre transfini~$1$. 
\citet{Dwork-1960} a obtenu la même conclusion si cette
série définit d'une part une fonction
méromorphe dans un disque complexe de rayon~$R_\infty$,
d'autre part une fonction méromorphe dans un disque $p$-adique de rayon~$R_p$
et si le produit des deux rayons $R_pR_\infty$ est~$>1$,
et \citet{Bertrandias-1963} a établi une  généralisation commune
des résultats de Pólya et Dwork.

Les preuves de ces théorèmes reposent sur le critère
de rationalité via les déterminants de Hankel, et sur des majorations
analytiques de ces déterminants telles que la formule
du produit entraîne leur nullité.

Ces énoncés ont été généralisés par Robinson, puis \citet{Cantor-1980},
lorsqu'on  considère plusieurs séries formelles basées
en plusieurs points distincts de la droite projective.
On exige alors que ces séries formelles
définissent des fonctions méromorphes complexes
ou $p$-adiques sur des domaines assez grands, ainsi qu'une propriété
d'intégralité adéquate, et la conclusion est qu'elles proviennent d'une
même fraction rationnelle.

Plus récemment, \citet{BostChambert-Loir-2009} ont étendu
le théorème de Borel--Dwork--Pólya--Bertrandias
en remplaçant la droite projective par une courbe algébrique de
genre arbitraire. Même dans le cas de la droite projective,
leur démonstration est fondamentalement différente.
Une première partie établit que la série formelle est
le développement de Taylor d'une fonction algébrique ;
elle utilise des techniques d'approximation diophantienne,
selon une méthode inaugurée par~\cite{ChudnovskyChudnovsky-1985},
puis généralisée par~\citet{Andre-2004a} (l'article date de 1997)
et dans un contexte non linéaire par \citet{Bost-2001}.
Une seconde partie en déduit la rationalité ; l'outil essentiel
est le théorème de l'indice de Hodge en théorie d'Arakelov
pour des fibrés en droites munis de métriques (hermitiennes ou $p$-adiques)
fournies par la théorie du potentiel.

Le but de cet article est d'établir
un théorème de rationalité (théorème~\ref{theo.rationnel})
du genre de ceux évoqués ci-dessus
qui étende simultanément le théorème de \cite{Cantor-1980}
et celui de \cite{BostChambert-Loir-2009}.

Nous reprenons l'approche de ce dernier article en
démontrant d'abord un énoncé d'algébricité  
(théorème~\ref{theo.algebrique}) pour des familles
finies de courbes formelles tracées sur une variété algébrique.
Le paragraphe~\ref{s-enonce} est consacré à mettre
en place l'énoncé de ce théorème qui fait intervenir
la \emph{valeur d'un jeu} défini à l'aide des théories du potentiel
complexes et non archimédiens. Dans l'espoir que ce critère
puisse être utilisé lorsque l'on met en jeu une  infinité de places
du corps de base, les gains de ce jeu sont à valeurs
dans~$\R\cup\{+\infty\}$.

Les définitions posées, nous prouvons ensuite ce théorème d'algébricité
au paragraphe~\ref{s-algebrisation},
en reprenant la méthode d'approximation diophantienne,
dans le langage de la méthode des pentes de~\cite{Bost-2001}.
Comme il y a plusieurs points, la méthode d'approximation diophantienne
doit dériver successivement par rapport à chacun d'entre eux,
comme dans la preuve du théorème de Schneider--Lang.
Dans ce cas, il faut également reprendre une innovation mise en place
par \citet{Herblot-2011} et dériver à des vitesses différentes
en chacun des points. 

Dans un dernier paragraphe, nous énonçons le 
le théorème de rationalité (théorème~\ref{theo.rationnel})
et expliquons comment
sa preuve se déduit du théorème d'algébricité et 
du théorème de l'index de Hodge en géométrie d'Arakelov.
Les arguments sont ici très proches de ceux de~\cite{BostChambert-Loir-2009}.

 
J'avais conçu le cœur de ces arguments il y a une dizaine
d'années et les avais mentionnés de-ci, de-là, sans pour autant
prendre le temps d'une rédaction en forme. 
Robert Rumely m'avait informé 
de la démonstration d'un tel théorème par \citet{Walters-2012}.
En septembre 2020, Vesselin Dimitrov m'avait également 
suggéré la possibilité d'un tel théorème.
Au moins si l'on se contente des places archimédiennes, 
on peut d'ailleurs le déduire des résultats de \cite{BostCharles-2026}.

Je remercie Mathilde Herblot, Robert Rumely, Vesselin Dimitrov, Jean-Benoît Bost et Guillaume Garrigos, de leur intérêt à des moments divers de ce travail.
Je remercie également les deux \emph{referees} pour leur lecture
attentive qui m'a conduit à corriger plusieurs insuffisances de la première
version et à ajouter quelques détails.

Camille Noûs est un chercheur français allégorique et polymathe, 
au genre indéfini, né en 2020 à l'initiative 
du groupe de défense de la recherche RogueESR.
Sa présence comme co-autrice de ce texte a pour but d'une part 
de protester contre les politiques court-termistes en matière de financement
de la recherche publique, en particulier leur dépendance toujours
plus grande à des appels à projet au sort incertain,
et d'autre part de reconnaître que ce texte s'inscrit
dans une longue chaîne de travaux de recherche auxquels je ne fais
qu'ajouter quelques commentaires.

Les quelques articles que j'ai pu
écrire avec Yuri Tschinkel depuis presque trente ans
m'ont été une grande source de stimulation intellectuelle.
C'est une grande joie que de publier cet article
dans un volume en son honneur.

\section{Théories du potentiel, jeux, algébricité}\label{s-enonce}

\subsection{Courbes A-analytiques}

Soit $F$ un corps de nombres et soit $R$ son anneau d'entiers.
On identifie l'ensemble des places de~$F$ à un système de représentants
vérifiant la formule du produit: pour toute place~$v$ de~$F$,
$\abs{\,\cdot\,}_v$ est une valeur absolue non triviale sur~$F$, dite valeur
absolue $v$-adique, et l'on suppose que l'on a $\prod_v \abs a_v=1$
pour tout $a\in F^\times$. 

Soit $X$ une variété algébrique quasi-projective sur~$F$.
Soit $(P_i)_{i\in I}$ une famille finie de points rationnels
de~$X$ et,
pour tout~$i$, soit $\widehat C_i\hookrightarrow \widehat X_{P_i}$ 
une courbe formelle lisse qui est \emph{A-analytique.}

Cette notion d'A-analyticité a été introduite par~\citet[déf.~3.7]{BostChambert-Loir-2009}
et requiert la conjonction de deux conditions:
\begin{enumerate}
\item Pour toute place~$v$ de~$F$, la courbe formelle $\widehat C_i$
se prolonge en un disque analytique sur~$F_v$;
\item 
Pour tout indice~$i$,
le produit des \emph{tailles} $v$-adiques $S_v(\widehat C_i)$
de~$\widehat C_i$,
quand $v$ parcourt l'ensemble des places finies de~$F$, 
est strictement positif.
\end{enumerate}
Nous renvoyons à \citep[\S3.1.1]{Bost-2001}
ou à \citep[\S3.A]{BostChambert-Loir-2009} pour
la définition précise de ces tailles.
Dit grossièrement, 
étant donné un plongement de~$X$
dans un schéma projectif et lisse~$\mathscr X$
sur l'anneau~$R$ des entiers de~$F$,
ces tailles mesurent le rayon d'un disque analytique
sur le corps~$F_v$  qui « prolonge » $\widehat C_i$.
Elles sont invariantes par immersion et par morphisme étale.
Si leur valeur précise dépend du choix de ce plongement,
la condition de stricte positivité de leur produit n'en dépend pas.

Donnons néanmoins trois exemples :
\begin{enumerate}
\item Si $\widehat C_i$ est un germe de courbe analytique $v$-adique,
sa taille $v$-adique est strictement positive.
\item Si $\widehat C_{i}$ est la fibre générique d'un
sous-schéma formel lisse de dimension~$1$ de~$\mathscr X_{R_v}$,
sa taille $v$-adique vaut~$1$; en particulier,
les localisations formelles de courbes algébriques en un
point lisse sont A-analytiques.
\item Si $\widehat C_{i}$ est une feuille formelle d'un feuilletage
de dimension~$1$ d'un modèle lisse $\mathscr X_{R_v}$,
sa taille $v$-adique est minorée par $\abs p_v^{1/(p-1)}$, où
$p$ est la caractéristique résiduelle de~$v$, 
et par $\abs p_v^{1/p(p-1)}$ si le feuilletage résiduel
est stable par $p$-courbures. En particulier, 
si $\widehat C_i$ est une feuille formelle d'un feuilletage qui 
vérifie la condition de Grothendieck--Katz,
elle est A-analytique.
\end{enumerate}

\subsection{Normes de référence}
Comme $\mathscr X$ est propre sur~$R$, les points rationnels $P_i\in X(F)$
se prolongent de manière unique 
en des sections $\mathscr P_i\colon \Spec(R)\to\mathscr X$.

Le modèle entier~$\mathscr X$ induit en outre, pour tout $i\in I$,
un sous-$R$-module de type fini de la $F$-droite~$T\widehat C_i$:
son dual est le sous-$R$-module de $T^*\widehat C_i$
engendré par $\mathscr P_i^*\Omega^1_{\mathscr X/R}$.

Pour toute place non archimédienne~$v$,
on en déduit un sous $R_v$-module de type fini de la $F_v$-droite~$T\widehat C_i\otimes_F F_v$ puis une norme $v$-adique sur cet espace vectoriel.

Lorsque $v$ parcourt les places archimédiennes de~$F$,
on se donne également des normes « $v$-adiques »
sur $T\widehat C_i\otimes_F F_v$,
invariantes par la conjugaison complexe.

On se donne enfin pour tout~$i$ une base $\xi_i$ de $T\widehat C_i$.

\subsection{Théories du potentiel}
Pour toute place~$v$ de~$F$, on suppose donnée
une « courbe analytique »~$\Omega_v$ sur~$F_v$
munie d'une famille finie de points $(z_i)$, deux à deux disjoints,
et un morphisme $\phi_v\colon\Omega_v\to X_v$
tel que, pour tout~$i$, $\phi_v$ définisse un isomorphisme
du complété formel de~$\Omega_v$ en~$z_i$ sur la courbe formelle
$\widehat C_{i,F_v}$. Plus précisément:
\begin{itemize}
\item Si $v$ est une place complexe, $\Omega_v$ est une surface
de Riemann,
et $\phi_v\colon \Omega_v\to X(\C_v)$ est une application holomorphe;
\item Si $v$ est une place réelle, $\Omega_v$ est une surface
de Riemann munie d'une involution
anti-holomorphe compatible avec $\phi_v$ et l'involution
anti-holomorphe canonique de $X(\C_v)$;
\item Si $v$ est non archimédienne, $\Omega_v$
est une courbe $F_v$ analytique lisse (sans bord),
au sens de~\cite{Berkovich-1990},
et $\phi_v\colon\Omega_v\to X_{F_v}^\an$ est un morphisme
d'espaces analytiques.
\end{itemize}

Par exemple, on peut prendre pour courbe~$\Omega_v$ 
la réunion disjointe, pour $i\in I$,
d'un disque ouvert « centré » en~$P_i$ et de rayon la taille $v$-adique
$S_v(\widehat C_i)$, 
de sorte que $\norm{\phi_{v,*}(\partial t/\partial t_i)}=1$
pour le paramètre local canonique de ces disques.

On peut en fait toujours supposer que~$\Omega_v$ 
contient ces disques.

Nous allons faire usage de  la théorie du potentiel sur une courbe analytique. 
Dans le cas complexe, nous renvoyons à l'article de~\cite{Bost-1999};
la théorie non archimédienne est initialement due à~\cite{Rumely-1989},
mais nous utiliserons la version de~\citet{Thuillier-2005}
dans le cadre des espaces de Berkovich.

Pour simplifier l'exposition, on suppose en fait
que $\Omega_v$ est l'intérieur d'une courbe analytique compacte~$\Omega'_v$
(au sens des surfaces de Riemann à bord si $v$ est
archimédienne, au sens des courbes de Berkovich si $v$ est non archimédienne)
et que $\Omega'_v\setminus\Omega_v$ n'est pas polaire.
Dans le cas non archimédien, il ne sera pas restrictif
de supposer que~$\Omega_v$ est le complémentaire d'un domaine
affinoïde non vide d'une courbe analytique propre.

Pour tout~$i$,
la \emph{fonction de Green}  $g_{\Omega_v,z_i}\colon \Omega'_v\setminus\{z_i\}\to\R$ 
est caractérisée par les propriétés suivantes:
\begin{enumerate}
\item Elle est continue sur~$\Omega'_v$, 
nulle sur $\Omega'_v\setminus \Omega_v$;
\item Elle est harmonique sur~$\Omega_v$;
\item Si $t_i$ est un paramètre local au voisinage de~$z_i$, 
alors $g_{\Omega_v,z_i}(z)+\log( \abs{t_i(z)} )$ a une limite en~$z_i$.
\end{enumerate}
Pour tout~$i$,
cette fonction de Green munit le fibré en droites
$\mathscr O_{\Omega_v}(z_i)$ sur~$\Omega'_v$
d'une métrique $v$-adique par la formule
\[  \log (\norm{\mathbf 1_{z_i} (z)}_v^{-1}) = g_{\Omega_v,z_i}(z), \]
Via la formule d'adjonction, la droite~$T\widehat C_i$ est alors
munie 
d'une \emph{norme capacitaire} $v$-adique, caractérisée par la formule:
\[ - \log(\norm{\phi_{v,*} (\partial/\partial t_i)}_v )
 = \lim_{z\to z_i} \left(g_{\Omega_v,z_i}(z)+\log( \abs{t_i(z)}_v) \right). \]
Le principe du maximum entraîne 
que $g_{\Omega_v,z_i}(z)$ décroît lorsque~$\Omega_v$ 
est remplacé par un ouvert contenant~$z_i$ et~$z$,
et que cette norme croît lorsque~$\Omega_v$ 
est remplacé par un ouvert plus petit.
Nous renvoyons à~\cite[3.1.4]{Bost-1999},
\cite[remarque~5.14, b)]{BostChambert-Loir-2009}
et \cite[\S3.6]{Thuillier-2005}
pour plus de détails.

Comme on a supposé que $\Omega_v$ contient le disque ouvert
de centre~$z_i$ et de rayon~$S_v(\widehat C_i)$,
on a les inégalités, pour tout $i\in I$,
\begin{equation}
g_{\Omega_v,z_i}(z) \geq \log (S_v(\widehat C_i)/\abs{t_i(z)}_v))
\end{equation}
et 
\begin{equation}
 -\log (\norm{\phi_{v,*}(\partial/\partial t_i)}_v) \geq 
\log (S_v(\widehat C_i)). \end{equation}

Pour toute place~$v$ de~$F$,
considérons alors la matrice~$\mathmat G^v$ de type~$I\times I$
dont le coefficient $(i,j)$ est donné par
\[  \mathmat G^v_{i,j}= g_{\Omega_v,z_i}(z_j) \]
si $j\neq i$,
et par
\[ 
\mathmat G^v_{i,i} = \lim_{z\to z_i} \left(g_{\Omega_v,z_i}(z) + \log (\abs{t_i(z)}_v) \right) \]
si $j=i$, où $t_i$ est un paramètre local sur~$\Omega_v$ en~$z_i$
tel que $(\phi_v)_*(\partial/\partial t_i)=\xi_i$.

Cette matrice ne dépend du choix des coordonnées locales~$t_i$
que par le choix des vecteurs~$\xi_i\in T\widehat C_i$.
Précisément, si $\xi_i$ est remplacé par~$\xi'_i=a_i\xi_i$,
on aura $\log (\abs{t'_i(z)})-\log(\abs{t_i(z)})\to -\log(\abs{a_i}_v)$
lorsque $z\to z_i$, de sorte que $\mathmat G^v$ est remplacée
par la matrice $\mathmat G^{\prime v}$ donnée par 
$\mathmat G^{\prime v}_{i,i} = \mathmat G^v_{i,i}-\log (\abs{a_i}_v)$
et $\mathmat G^{\prime v}_{i,j} = \mathmat G^v_{i,j}$ si $j\neq i$.

La matrice~$\mathmat G^{v}$ est symétrique: voir
\cite[p.~80]{Rumely-1989} dans le cas archimédien et
\cite[remarque 3.4.15]{Thuillier-2005} dans le cas non archimédien.

\subsection{Capacités globales} 
%

Comme les courbes~$\widehat C_i$ sont A-analytiques,
la somme~$\mathmat G$ des matrices~$(\mathmat G^v)$,
lorsque $v$ parcourt  l'ensemble des places de~$F$,
est une matrice de type~$I\times I$ à coefficients
dans $\R\cup\{+\infty\}$, symétrique,
et positive ou nulle (éventuellement infinie) hors de la diagonale.

Notons qu'au contraire des matrices locales~$\mathmat G^v$,
la matrice globale~$\mathmat G$  
ne dépend plus du choix des vecteurs tangents~$\xi_i$:
lorsque $\xi_i$ est changé en~$\xi'_i=a_i\xi_i$,
le coefficient~$\mathmat G_{ii}$ est changé en
\[ \mathmat G'_{ii} = \mathmat G_{ii}- \sum_v \log (\abs{a_i}_v)
 = \mathmat G_{ii} \]
en vertu de la formule du produit.

\subsection{Théorie des jeux}
Soit 
\[  V_{\mathmat G} = \Sup_{x\in \Delta_I} \Inf_{y\in \Delta_I} \langle x, \mathmat Gy \rangle \]
la « valeur du jeu » donné par cette matrice,
où $\Delta_I$ est le simplexe des $x\in\R^+_I$ tels que $\sum x_i=1$.
En théorie des jeux, les lignes  et les colonnes de la matrice~$\mathmat G$
représentent les coups possibles pour les deux joueurs,
Isaac et Jacob, 
et le nombre réel~$\mathmat G_{i,j}$ (éventuellement infini si $i\neq j$)
situé sur la ligne~$i$ et la colonne~$j$
représente le gain d'Isaac
lorsqu'il joue le coup~$i$ et que 
Jacob
répond par le coup~$j$. Un point $x\in \Delta_I$ est classiquement 
interprété comme une « stratégie mixte » revenant à jouer
le coup~$i$ avec probabilité~$x_i$. Dans la formule
précédente pour~$V_{\mathmat G}$, l'expression $\langle x, \mathmat Gy\rangle$
représente le gain moyen 
d'Isaac
lorsqu'il adopte la stratégie mixte~$x$
et que 
Jacob
adopte la stratégie mixte~$y$;
la borne inférieure sur~$y$ 
suivie de la borne supérieure sur~$x$ signifie que
la stratégie mixte de Jacob est de 
minimiser le gain d'Isaac,
et qu'Isaac cherche à le maximiser.
Lorsque $\mathmat G$ est à valeurs finies, 
le théorème du minimax de \cite{vonNeumann-1928}
affirme que 
la valeur~$V_{\mathmat G}$ est également donnée
par la formule (a priori plus grande)
\[  V_{\mathmat G} =  \Inf_{y\in \Delta_I} \Sup_{x\in \Delta_I}\, \langle x, \mathmat Gy \rangle. \]
Cette commutativité traduit l'\emph{équilibre} du jeu:
$V_{\mathmat G}$ est effectivement le gain que peut espérer 
Isaac
face à toute stratégie mixte de Jacob.

\begin{theo}\label{theo.algebrique}
Supposons  $V_{\mathmat G}>0$. Alors les courbes formelles~$\widehat C_i$
sont algébriques: il existe une plus petite
courbe algébrique $Y\subset X$
telle que pour tout $i$, $P_i\in Y(F)$ 
et $\widehat C_i \subset \widehat Y_{P_i}$.
\end{theo}

\begin{remas}
\phantomsection
\label{remas.algebrique}
\begin{enumerate}
\item
Dans le cas particulier où il n'y a qu'une courbe formelle~$C_i$,
l'ensemble~$I$ est réduit à l'élément~$i$,
la matrice~${\mathmat G}$ se réduit à 
la somme des coefficients $\mathmat G^v_{ii}$
qui fournissent des normes $v$-adiques
sur $T\widehat C_i$, 
analogues aux métriques $v$-adiques capacitaires
de~\cite{BostChambert-Loir-2009},
de sorte que $\mathmat G= \hdeg (T\widehat C_i)$,
puis $V_{\mathmat G}=\hdeg (T\widehat C_i)$.
On retrouve alors le théorème~6.1 de \citet{BostChambert-Loir-2009}.

En général, si $i\in I$,
en prenant pour $y$ l'unique élément de~$\Delta_I$
tel que $y_i=1$ et en utilisant
que les coefficients non diagonaux de la matrice~${\mathmat G}$  sont positifs
ou nuls, on obtient $V_{\mathmat G}\geq \hdeg (T\widehat C_i)$.
Autrement dit, dans le cadre du théorème~\ref{theo.algebrique}
il est a priori plus efficace de prendre en compte
que l'on a plusieurs points à disposition.

\item
La matrice~${\mathmat G}$ décroît lorsque les~$\Omega_v$ 
sont remplacés par des ouverts. La
définition de~$V_{\mathmat G}$  prouve qu'elle 
décroît également. Le théorème obtenu est d'autant plus
« facile » à appliquer que les ouverts~$\Omega_v$ sont grands.

\item
Dans certains cas, on peut affirmer que la courbe~$Y$ est intègre:
c'est en particulier le cas s'il existe une place~$v$
telle que la courbe~$\Omega_v$ soit connexe,
Considérons en effet une composante irréductible~$Z$ de~$Y$
et prouvons que $Z=Y$. Pour cela, il suffit de démontrer
que pour tout couple $(i,j)$ dans~$I$ tel que $\widehat C_i\subset Z$,
on a $\widehat C_j\subset Z$.
Soit $f$ une fonction algébrique sur (un ouvert
affine de~$X$ contenant~$P_j$) qui est identiquement nulle sur~$Z$ et
démontrons qu'elle s'annule sur~$\widehat C_j$.
Par hypothèse, la fonction analytique~$\phi_v^*f$ sur~$\Omega_v$
s'annule à  l'ordre infini en~$z_i$, 
donc est identiquement nulle sur~$\Omega_v$, puisque $\Omega_v$
est supposée connexe. Par suite, $\phi_v^*f$ s'annule
sur~$\widehat C_j$, et donc $f$ aussi.

Si $v$ est une place de~$F$ et si $i$ et~$j$ sont deux indices distincts
tels que $\mathmat G^v_{i,j}>0$, le principe du maximum
entraîne que~$z_i$ et~$z_j$ sont contenus dans une même
composante connexe de~$\Omega_v$ ; 
la réciproque est d'ailleurs vraie si $v$ est archimédienne.
Ainsi, si $\mathmat G^v_{i,j}>0$ pour tous~$i$
et~$j$ tels que $i\neq j$, la courbe~$Y$ est irréductible.

\item
Supposons maintenant que la matrice globale~$\mathmat G$ 
est irréductible au sens de la théorie des chaînes de Markov
et démontrons que la courbe~$Y$ est irréductible.
Considérons de nouveau une composante irréductible~$Z$ de~$Y$
et démontrons que $Z=Y$.
Par définition d'une matrice irréductible, il suffit de démontrer
que pour tout couple~$(i,j)$ dans~$I$ tel que $i\neq j$,
$\mathmat G_{i,j}>0$ et $\widehat C_i\subset Z$, on a 
$\widehat C_j\subset Z$.
Par hypothèse, il existe une place~$v$ telle que $g_{\Omega_v,z_i}(z_j)>0$,
et les points~$z_i,z_j$ appartiennent à la même composante connexe
de~$\Omega_v$. On en déduit comme précédemment que $\widehat C_j\subset Z$.
\end{enumerate}
\end{remas}

\section{Algébrisation}
\label{s-algebrisation}
 
 
La démonstration du théorème~\ref{theo.algebrique} utilise une technique
d'approximation diophantienne sous la forme de la méthode
des pentes inventée par~\cite{Bost-1996}, mise
en œuvre dans cette direction par~\cite{Bost-2001}
et~\cite{BostChambert-Loir-2009}.
On utilise
également un raffinement inauguré par~\cite{Herblot-2011}
dans le contexte de la méthode des pentes
qui permet de ne pas dériver à la même vitesse en chaque point.

\subsection{Réductions}
Il n'est pas restrictif de supposer que $X$ soit projectif.
Soit $Y$ le plus petit sous-schéma fermé de~$X$ tel que pour tout~$i$,
on ait $P_i\in Y(F)$ et $\widehat C_{i}\subset\widehat Y_{P_i}$.
Quitte à remplacer~$X$ par~$Y$, on peut supposer que $Y=X$. 
Il s'agit alors de démontrer que $\dim(X)=1$. 

\begin{lemm}\label{lemm.minimax}
Soit $V'_{\mathmat G}$ un nombre réel tel que $V'_{\mathmat G}<V_{\mathmat G}$.
Il existe un élément $a\in\Delta_I$ dont les coordonnées
sont des nombres rationnels strictement positifs
tel que $\sum_{i} a_i \mathmat G_{ij}>V'_{\mathmat G}$ pour tout $j\in I$.
\end{lemm}
\begin{proof}
Par définition d'une borne supérieure, 
il existe $x\in\Delta_I$ 
tel que $\langle x,{\mathmat G}y\rangle > V'_{\mathmat G}$
pour tout~$y\in\Delta_I$, soit, de manière équivalente 
$\sum_i x_i {\mathmat G}_{ij} > V'_{\mathmat G}$ pour tout $j\in I$.
Les fonctions $x\mapsto \sum_i x_i \mathmat G_{i,j}$
sont semi-continues inférieurement; par suite, 
tout vecteur~$a\in\Delta_I$ qui est assez proche de~$x$ 
vérifiera encore $\sum a_i \mathmat G_{ij}\geq V'_{\mathmat G}$
pour tout~$j$.
Il en existe en particulier dont les coordonnées sont rationnelles
et strictement positives.
\end{proof}

\subsection{}
Rappelons que nous avons fixé un plongement de~$X$ dans
un $R$-schéma projectif et lisse~$\mathscr X$,
munissons-le d'un fibré en droites hermitien ample~$\mathscr L$;
notons~$L$ la restriction de~$\mathscr L$ à~$X$.
Notons~$\mathscr X'$ l'adhérence de~$X$ dans~$\mathscr X$.

Pour tout entier naturel~$N$, l'espace des
sections $\mathscr E_N=\Gamma(\mathscr X',\mathscr L^{\otimes N})$
est un $R$-module projectif de type fini, 
de fibre générique $E_N=\Gamma(X,L^{\otimes N})$.

À toute place archimédienne~$\sigma$ de~$R$,
l'espace~$E_N$ possède une norme naturelle
(« du sup »), notée~$\norm{\,\cdot\,}_{\sigma,\infty}$,
déduite de l'évaluation ponctuelle des sections,
et nous le munirons de l'unique norme hermitienne~$\norm{\,\cdot\,}_\sigma$
dont la boule unité a volume minimal parmi les normes qui majorent
la norme~$\norm{\,\cdot\,}_{\sigma,\infty}$
(« norme de John », voir~\cite{Gaudron-2008,Chen-2012a}).
Cela fait de~$\mathscr E_N$ un $R$-fibré vectoriel hermitien~$\overline{\mathscr E_N}$.

Si $n=\dim(X)$,
on rappelle le comportement asymptotique, quand $N\to+\infty$,
du rang et du degré arithmétique de ce fibré vectoriel hermitien:
comme $L$ est ample, le théorème de Hilbert--Samuel s'écrit:
\begin{equation}\label{eq.HS}
\operatorname{rang}(\mathscr E_N) \sim \deg_L(X)\, N^n.
\end{equation}
Le théorème de Hilbert--Samuel arithmétique 
\citep{GilletSoule-1988, AbbesBouche-1995, Zhang1995b}
calcule la limite de la suite $(N^{-n-1}\hdeg(\overline{\mathscr E_N}))$; 
on en déduit l'existence d'un nombre réel~$c$ tel que 
\begin{equation}\label{eq.aHS}
\hdeg(\overline{\mathscr E_N}) \geq -c N^{n+1}.
\end{equation} 
Signalons aussi que cette minoration se démontre élémentairement 
\citep[proposition~4.4]{Bost-2001}. 

\subsection{Évaluations}
Nous allons définir une filtration décroissante et exhaustive $(E^k_N)_k$
de~$E_N$.
Dans la méthode des pentes avec un seul point d'évaluation~$P$,
le $k$-ième cran de la filtration~$E^k_N$ correspond aux
éléments~$s\in E_N$ dont la restriction à~$\widehat C_P$
s'annule au moins à l'ordre~$k$.

Disposant d'une famille de points~$(P_i)_{i\in I}$,
nous allons nous donner une suite $(i_k)_{k\geq1}$ dans~$I$,
poser $E^0_N=E_N$ et, pour $k\geq 1$, 
définir~$E^{k}_N$ comme l'ensemble des éléments de~$E^{k-1}_N$
dont la restriction à~$\widehat C_{{i_k}}$
s'annule une fois de plus que ce qui était imposé
au cran précédent. 
Définissons des suites $(\omega_i(k))_k$ par
\[ \omega_i(k)=\Card(\{ j\in\{1,\dots,k\}\sozat i_j=i \})
= \sum_{1\leq j\leq k}  \delta_{i_j,i}, \]
où $\delta$ désigne l'indice de Kronecker;
on a donc 
$\omega_{i_k}(k)=\omega_{i_k}(k-1)+1$, et $\omega_i(k)=\omega_i(k-1)$
si $i\neq i_k$ ; 
on a en particulier, pour tout $k\in\N$, la relation
\[ \sum_{i\in I} \omega_i(k) = k. \]
Enfin,  pour tout entier naturel~$k$,
le sous-module~$E^{k}_N$ de~$E_N$ est l'ensemble
des éléments~$s$ de~$E_N$
tels que pour tout~$i\in I$, la restriction de~$s$ à $\widehat C_{i}$
s'annule au moins à l'ordre~$\omega_i(k)$.

Par définition, pour $k\geq 1$, le sous-espace $E^{k}_N$
est le noyau du \emph{morphisme d'évaluation}
\[ \phi^k_N\colon E^{k-1}_N\to T^*\widehat C_{{i_{k}}}^{\otimes \omega_{i_{k}}(k-1)}\otimes L^{\otimes N}(P_{i_{k}}). \]

Définissons maintenant la suite~$(i_k)$.

Soit $V'_{\mathmat G}$ un nombre réel tel que $0<V'_{\mathmat G}<V_{\mathmat G}$.
D'après le lemme~\ref{lemm.minimax}, il existe
un vecteur $a\in\Delta_I$ 
à coordonnées strictement positives et rationnelles
tel que $({\mathmat G}a)_i>V'_{\mathmat G}$ pour tout~$i$.

On définit alors la suite $(i_k)_k$ par récurrence comme suit:
si $i_{1},\dots,i_{k}$ sont choisis,
on prend  pour $i_{k+1}$ l'un des indices~$i$
pour lequel $\omega_i(k)-k a_i$ est minimal.

La construction garantit que les suites $(\omega_i(k)-ka_i)_k$ sont bornées.
Plus précisément :
\begin{lemm}\label{lemm.omegai(k)}
Pour tout $k\in\N$ et tout $i\in I$, on a 
\[  1-\Card(I)\leq \omega_i(k)-ka_i \leq 1 .\]
\end{lemm}
\begin{proof}
On raisonne par récurrence sur~$k$. Ces inégalités sont évidentes
lorsque $k=0$, car le terme central est alors nul; 
supposons $k>0$ et que ces inégalités soient satisfaites au rang~$k-1$.
Puisque $\sum_{i\in I}(\omega_i(k-1)-(k-1)a_i)=0$, 
le choix de~$i_{k}$ entraîne que $\omega_{i_{k}}(k-1)-(k-1)a_{i_{k}}\leq 0$, 
d'où \[ \omega_{i_{k}}(k)-ka_{i_{k}}
=(\omega_{i_{k}}(k-1)-(k-1) a_{i_{k}})+(1-a_{i_{k}})\leq 1-a_{i_k} \leq 1.\]
Pour tout autre indice~$i$, on a 
\[ \omega_i(k)-ka_i=(\omega_i(k-1)-(k-1) a_i ) -a_i
\leq  \omega_i(k-1)-(k-1) a_i\leq1,\]
compte tenu de l'hypothèse de récurrence.
Puisque la somme, pour $i\in I$, des $\omega_i(k)-ka_i$ est nulle
et que chacun est majoré par~$1$, ils sont tous minorés par
$1-\Card(I)$.
\end{proof}

\begin{coro}\label{coro.exhaustive}
La filtration~$(E_N^k)_k$ de~$E_N$ est finie et exhaustive.
\end{coro}
\begin{proof}
D'après les inégalités du lemme précédent, les suites~$(\omega_i(k))_k$
tendent vers~$+\infty$, car $a_i>0$ pour tout~$i$. Ainsi, 
la restriction à~$\widehat  C_i$ 
d'un élément~$s$ de~$\bigcap_k E_N^k$ est nulle. 
Par définition de~$X$, un tel élément est nul.
\end{proof}

\begin{coro}\label{coro.minor}
Il existe un nombre réel~$c$ tel que pour tout $k\in\N$
et tout $j\in I$, on ait
\[ \sum_{i\in I} \omega_i(k) \mathmat G_{ij} \geq k V'_{\mathmat G} -c. \]
\end{coro}
\begin{proof}
On écrit en effet
\[
 \sum_{i\in I} \omega_i(k) \mathmat G_{ij}
 = k \sum_{i\in I} a_i \mathmat G_{ij}
+ \sum_{i\in I} (\omega_i(k)-k a_i) \omega_i(k) .\]
D'après le choix des~$a_i$, le premier terme est minoré par $kV'_{\mathmat G}$.
D'après le lemme précédent, le second est borné lorsque $k$ varie.
\end{proof}

\subsection{}
Puisque $\mathscr X$ est propre sur~$\Spec(R)$, 
pour tout~$i$, le point $P_i\in X(F)$ donne lieu à une
section $\mathscr P_i\colon\Spec(R)\to\mathscr X$.
L'image de $\mathscr P_i^*\Omega^1_{\mathscr X/R}$
dans~$T^*\widehat C_{i}$ est
un sous-$R$-module projectif de rang~$1$ de~$T^*\widehat C_{i}$
que l'on note $(T^*\widehat C_i)_{\mathscr X}$.
On le munit également de métriques hermitiennes aux places
archimédiennes de~$F$.

Dans ces conditions, pour tout $i\in I$ et tout $a\in\N$,
le $R$-module projectif de rang~$1$
\[ (T^*\widehat C_{i})_{\mathscr X}^{\otimes a}\otimes \mathscr P_i^* \mathscr L^{\otimes N}, \]
est muni de métriques hermitiennes naturelles.
Son degré d'Arakelov est égal à
\[ - a \hdeg(T\widehat C_{i})_{\mathscr X} + N \hdeg(\mathscr P_{i}^*\overline{\mathscr L}). \]

\subsection{}
Compte tenu des choix faits, le morphisme d'évaluation~$\phi^k_N$
possède, pour toute place~$v$ de~$F$, une semi-norme $v$-adique
(réelle ou complexe si $v$ est archimédienne).
Nous noterons $h_v(\phi^k_N)$
le logarithme 
de cette semi-norme; c'est un élément de $\R\cup\{-\infty\}$,
négatif ou nul pour presque toute place~$v$.
Posons aussi 
\[ h(\phi^k_N) = \sum_v h_v(\phi^k_N). \]
C'est un élément de $\R\cup\{-\infty\}$.

 
\begin{prop}\label{prop.major-locale}
Soit $v$ une place de~$F$. 
\begin{enumerate}
\item
Si $v$ est non archimédienne, on a 
\[
h_v(\phi^k_N) \leq -\omega_{i_k}(k) \log(S_{\mathscr X,v}(\widehat C_{i_{k+1}})).
\]
\item
Pour tout nombre réel~$\eps>0$, il existe un nombre réel~$c_v$
tel que l'on ait 
\[
 h_v(\phi^k_N) \leq 
- \sum_{i} \omega_i(k-1) \mathmat G^v_{i,i_k} 
+ N c_v + (k-1) \eps 
- \omega_{i_k}(k-1) \log (\norm{\mathrm d t_{i_k}}_v)
\]
pour tout~$N\in\N$ et tout~$k\in\N$.
\end{enumerate}
\end{prop}

\begin{proof}
La première majoration est une conséquence directe
du lemme~3.3 de~\citet{Bost-2001};
voir aussi~\citep[(6.3)]{BostChambert-Loir-2009}.

La seconde est une variante 
de \citep[prop.~3.6]{Bost-2004b} dans le cas archimédien,
et de \citep[prop.~5.15]{BostChambert-Loir-2009} dans le cas
ultramétrique,
mais sa démonstration requiert de nouveaux arguments.

Soit $\eps$ un nombre réel~$>0$.
Puisque les fonctions de Green~$g_{\Omega_v,z_i}$ tendent vers~$0$ à l'infini,
il existe dans~$\Omega_v$,
\begin{itemize}
\item
Si $v$ est archimédienne, une surface de Riemann compacte à bord~$W$ n'ayant aucune composante connexe sans bord;
\item
Si $v$ est non archimédienne, un domaine analytique compact~$W$ sans composante connexe propre;
\end{itemize}
en dehors duquel ces fonctions de Green sont toutes majorées par~$\eps$; 
en particulier, $W$ contient les points~$z_i$.
Dans le cas non archimédien, $W$ est alors affinoïde ;
comme le corps résiduel de~$F_v$ est fini, on déduit
de \citep[proposition~3.1]{vanderPut-1980};
que son groupe de Picard est de torsion et 
on choisit alors un entier~$d\geq1$ tel que le fibré
en droites $\phi_v^*L|_W^{\otimes d}$ soit trivial. 
Dans le cas archimédien, le groupe de Picard de~$W$ est trivial 
et on pose $d=1$.
(Comme aucune composante connexe de~$W$ n'est sans bord,
le double~$W'$ de~$W$ est une surface de Riemann connexe compacte ;
un fibré en droites holomorphe sur~$W$ fournit, par doublement, un fibré en droites holomorphe sur~$W'$ ; si $p'$ est un point de~$W'\setminus W$,
la surface de Riemann $W'\setminus\{p'\}$ n'est pas compacte,
et le théorème~30.3 de~\cite{Forster-1981} affirme que tout fibré en
droites qu'elle porte est trivial.)

Choisissons maintenant une trivialisation $\sigma$ de $\phi_v^*L|_W^{\otimes d}$.
Comme $W$ est compact, la norme de~$\sigma$ est bornée inférieurement
et supérieurement sur~$W$; 
posons $c_v=\frac1d \sup(\log (\norm{\sigma}))-\frac1d \inf(\log (\norm{\sigma}))$.

Après ces préparatifs, considérons maintenant $s\in E^{k-1}_N$. 
Pour tout $i\in I$, la restriction de~$s$ à~$\widehat C_i$
s'annule à l'ordre~$\omega_i(k-1)$ en~$P_i$, par définition de~$E^{k-1}_N$,
donc $\phi_v^*s$ s'annule à l'ordre~$\omega_{i}(k-1)$ en~$z_{i}$.
En particulier,
le développement formel de $\phi_v^*s$ au voisinage du
point~$z_{i_{k}}$ s'écrit 
\[ \phi_v^* s 
  = u  \phi_v^*t_{i_k}^{\omega_{i_{k}}(k-1)} + \cdots, \]
où $u\in L^N(P_{i_{k}})$,
et la définition du morphisme d'évaluation~$\phi^k_N$ entraîne
l'égalité
\[ \phi^k_N (s) =  (\mathrm d t_{i_k}) ^{\omega_{i_{k}}(k-1)} \otimes u 
 \]
dans $T^*\widehat C_i^{\otimes \omega_{i_k}(k)} \otimes L(P_{i_k})^{\otimes N}$.

Nous allons majorer~$\norm u$ à l'aide du principe du maximum.
Observons que $\phi_v^*s^d \otimes \sigma^{-N}$ 
est une fonction analytique sur~$W$ qui, pour tout~$i$,
s'annule à l'ordre~$d\,\omega_i(k-1)$ en~$z_i$.
La fonction 
\[ \log (\abs{\phi_v^*s^d\otimes \sigma^{-N}}) + d \sum_i \omega_i(k-1) g_{\Omega_v,z_i} \]
sur~$W$ est alors sous-harmonique, continue, et majorée par 
\begin{multline*}
 d\log (\norm s_{v,\infty})
 +  N \sup_{z\in W}\log(\norm{\sigma^{-1}(z)}) + d \sum_i \omega_i(k-1) \eps  \\
 = d\, \log (\norm s_{v,\infty})
 -  N \inf_{z\in W} \log(\norm{\sigma(z)})  + d\,(k-1)\eps 
\end{multline*}
au bord de~$W$.
D'après le principe du maximum, cette expression
majore sa valeur en~$z_{i_k}$, qui vaut, par définition, 
\begin{align*}
 d\, \log (\norm u) -  N \log (\norm {\sigma(z_{i_{k}})}) \hskip-5cm \\
& \qquad + d\,\omega_{i_{k}}(k-1) 
 \lim_{z\to z_{i_{k}}}  \left( g_{\Omega_v,z_{i_{k}}}(z) + \log(\abs{\phi_v^*t_{i_k}}(z))\right)  \\
& \qquad  
 + d\, \sum_{i\neq i_{k}} \omega_i(k-1) g_{\Omega_v,z_i}(z_{i_{k}})  \\
& =  d\,\log (\norm u)-N  \log (\norm{\sigma(z_{i_k})})
 + d \sum_i \omega_i(k-1) \mathmat G^v_{i,i_k}, 
\end{align*}
de sorte que 
\[ \log(\norm u) 
\leq N c_v
- \sum_{i}\omega_i(k-1) \mathmat G^v_{i,i_k} 
+ (k-1) \eps + \log(\norm{s}_{v,\infty})
.\]
\emph{In fine,} on a donc
\begin{align*}
h_v(\phi^k_N)  & \leq \log(\norm u)
+ \omega_{i_k}(k-1) \log(\norm{\mathrm d t_{i_k})})
- \log(\norm s_{v,\infty}) \\
&  \leq N c_v + (k-1) \eps 
- \omega_{i_k}(k-1) \log (\norm{\mathrm d t_{i_k})}) 
- \sum_i \omega_i(k-1) \mathmat G^v_{i,i_k}. 
\end{align*}
Cela conclut la preuve de la seconde
inégalité de la proposition~\ref{prop.major-locale}.
%
\end{proof}

\begin{coro}\label{coro.major-globale}
Il existe un nombre réel~$c'$ tel que,
pour tout nombre réel~$\eps>0$,
il existe un nombre réel~$c$ tel que l'on ait, pour tout $N\in\N$
et tout $k\in\N$, l'inégalité
\[ h(\phi^k_N) +\omega_{i_k}(k-1)\hdeg (T^*\widehat C_{i_k})
 \leq -kV'_{\mathmat G} + k c' \eps + N c  . 
\]
\end{coro}
\begin{proof}
Rappelons que nous avons choisi la famille~$(a_i)_{i\in I}$
de sorte que $\sum_i a_i \mathmat G_{ij} > V'_{\mathmat G}$ pour tout~$j\in I$.
Il existe donc un ensemble fini~$S$ de places de~$F$
contenant les places archimédiennes tel que
\[ \sum_{v\in S} \sum_i a_i \mathmat G^v_{ij} > V'_{\mathmat G} \]
et
\[ - \sum_{v\not\in S} \log (S_{\mathscr X,v}(\widehat C_j)) \leq\eps \]
pour tout $j\in I$.
On suppose également que $\norm{\mathrm dt_i}_v=1$
pour tout~$i\in I$ et toute place~$v\not\in S$.
En appliquant les estimées de la proposition~\ref{prop.major-locale},
suivant que $v\not\in S$, ou que $v\in S$ (en remplaçant~$\eps$
par~$\eps/\Card(S)$),
on en déduit qu'il existe une famille finie~$(c_v)_{v\in S}$
de nombres
réels telle que
\begin{align*}
 h(\phi^k_N) & = \sum_v h_v(\phi^k_N) \\
&  \leq - \omega_{i_k}(k-1)  \sum_{v\notin S} \log (S_{\mathscr X,v}(\widehat C_{i_{k+1}})) \\
& \qquad 
- \sum_{i\in I}  \omega_i(k-1) \sum_{v\in S} \mathmat G^{v}_{i,i_k}
- \omega_{i_k}(k-1) \sum_{v\in S} \log(\norm{\mathrm dt_{i_k}}_v \\
& \qquad\qquad + N \sum_{v\in S} c_v + (k-1) \eps .
\end{align*}
Dans cette expression, on commet une erreur bornée
en remplaçant $\omega_i(k-1)$ par $ka_i$ (lemme~\ref{lemm.omegai(k)}).
On obtient alors une majoration de la forme
\[ h(\phi^k_N) 
+ \omega_{i_k}(k-1) \hdeg ((T^*\widehat C_{i_k})_{\mathscr X})
\leq -k V'_{\mathmat G}+ O(k\eps)  + \mathrm O(N), \]
comme il fallait démontrer.
\end{proof}

\subsection{Inégalité de pentes}
L'\emph{inégalité de pentes} de~\cite[prop.~4.6]{Bost-1996}
associée à la filtration~$(E^k_N)$
de~$E_N$ et aux morphismes~$\phi^k_N$ s'écrit alors
\begin{multline}
  \hdeg(\overline{\mathscr E_N})
 \leq \sum_{k=1}^\infty  \rang(E^{k-1}_N/E^{k}_N) \times \\
\left(
 h(\phi^k_N)  +\omega_{i_{k}}(k-1) \hdeg((T^*\widehat C_{{i_{k}}})_{\mathscr X} )
+ N \hdeg(\mathscr P_{i_{k}}^*\overline{\mathscr L})\right). \end{multline}
Comme la filtration~$(E_N^k)_k$ est exhaustive (corollaire~\ref{coro.exhaustive}),
la série est en fait une somme finie.
En reportant l'inégalité du corollaire~\ref{coro.major-globale},
on obtient   donc
\begin{equation} 
\hdeg(\mathscr E_N) 
\leq \sum_{k=1}^\infty \rang(E^{k-1}_N/E^{k}_N) 
  (-k (V'_{\mathmat G}-c'\eps) + N c + N \hdeg(\mathscr P_{i_k}^*\overline{\mathscr L})).
\label{eq.pentes-2}\end{equation}

Choisissons~$\eps>0$ assez petit pour que $V'_{\mathmat G}-c'\eps >0$.
Compte-tenu des estimations~\eqref{eq.HS} et~\eqref{eq.aHS}, il existe
donc un nombre réel~$\kappa$  tel que
\begin{equation}
\label{eq.pentes-3}
\sum_{k=1}^\infty k \rang(E^{k-1}_N/E^{k}_N) \leq
\kappa N^{n+1} . 
\end{equation}
Nous allons en déduire que $n=1$, c'est-à-dire que $X$ est une courbe.

Par sommation d'Abel, on a en effet
\[ \sum_{k=1}^\infty k \rang(E^{k-1}_N/E^{k}_N)
= \sum_{k=0}^\infty  \rang(E^k_N). \]
D'autre part, puisque $\rang(E^{k-1}_N/E^{k}_N)\leq 1$,
on a 
\[ \rang(E^k_N) \geq \sup(\rang(E_N)-k,0), \]
d'où
\begin{align*}
 \sum_{k=1}^\infty k \rang(E^{k-1}_N/E^{k}_N)
& \geq \sum_{k=1}^{\rang(E_N)}  (\rang(E_N)-k)  \\
& = \frac12 \rang(E_N) (1+\rang(E_N)). \end{align*}
L'estimation~\eqref{eq.HS} et l'inégalité~\eqref{eq.pentes-3} 
impliquent alors une
majoration $N^{2n} = \mathrm O(N^{n+1})$, d'où $n=1$,
ce qui conclut la démonstration du théorème~\ref{theo.algebrique}.

\section{Rationalité}

\subsection{}
Le second théorème de cet article est un énoncé de rationalité
de fonctions formelles. 
On considère une courbe~$X$, lisse, géométriquement connexe et projective
sur le corps de nombres~$F$,
une famille finie $(P_i)_{i\in I}$ de points rationnels de~$X$ et,
pour tout $i\in I$, une fonction formelle~$f_i$ sur~$X$ le long de~$P_i$,
c'est-à-dire un élément de $\widehat{\mathscr O_{X,P_i}}$.

On suppose que pour toute place~$v$ de~$F$,
il existe un voisinage analytique~$\Omega_v$ de la réunion des~$P_i$
dans la courbe $F_v$-analytique~$X_v$ et
une fonction méromorphe~$f_v$ sur~$\Omega_v$
dont le développement de Taylor en~$P_i$  coïncide avec~$f_i$,
pour tout $i\in I$.

On fait en outre une hypothèse de nature adélique,
qu'il existe un modèle propre et lisse~$\mathscr X$ de~$X$
sur un sous-anneau de type fini~$R$ de~$F$ de corps des fractions~$F$
et des sections $\mathscr P_i\colon \Spec(R)\to \mathscr X$
prolongeant les points~$P_i$, de sorte que, d'une part chaque~$f_i$
définisse une fonction formelle de $\widehat{\mathscr X}_{\mathscr P_i}$,
et d'autre part pour toute place finie~$v$ dominant~$R$,
la courbe~$\Omega_v$ soit l'ouvert de~$X_v$ formé des points
ayant même image,  par l'application
de réduction $X_v \to \mathscr X_{k_v}$, que l'un des points~$P_i$.
(Ici, $k_v$ est le corps résiduel de~$R$ en~$v$.)

\subsection{}
Comme aux paragraphes précédents, les diverses théories du potentiel
fournissent des matrices~$\mathmat G^v$ de type~$I\times I$,
pour toute place~$v$, d'où on tire, de la même façon, leur somme~$\mathmat G$
et la valeur du jeu~$V_{\mathmat G}$.
Toutefois, ces matrices disposent maintenant d'une interprétation arithmétique 
plus claire: pour toute place~$v$ et tout indice~$i$, la fonction de Green
$g_{\Omega_v,P_i}$ définit une métrique $v$-adique 
sur le fibré en droites $\mathscr O_X(P_i)$ sur~$X$, et
$\mathmat G_{i,j}$ est le degré d'Arakelov 
de $\mathscr P_j^*\overline{\mathscr O_X(P_i)}$.
Ainsi, $\mathmat G$ est la matrice d'un « accouplement 
de hauteurs capacitaire » associé aux~$\Omega_v$ et aux points~$P_i$.

\begin{theo}\label{theo.rationnel}
Supposons $V_{\mathmat G}>0$. Il existe alors une fonction rationnelle~$f$
sur~$X$ de développement formel~$f_i$ en~$P_i$, pour tout $i\in I$.
\end{theo}

Dans le cas particulier où $X$ est la droite projective,
on retrouve la première partie du théorème~5.2.2 de~\cite{Cantor-1980}.
(La seconde partie de ce théorème affirme que la
condition $V_{\mathmat G}>0$ est nécessaire; nous n'avons
pas étudié sa généralisation au cas des courbes de genre arbitraire.)

Lorsque $X$ est la droite projective et où l'ensemble~$I$
est réduit à un élément, on obtient le classique
théorème de Pólya--Bertrandias, lui-même une élaboration
du théorème de Borel--Dwork.

Lorsque l'ensemble~$I$ est réduit à un élément, on retrouve
le théorème~7.9 de~\cite{BostChambert-Loir-2009}
dont nous reprenons le schéma de démonstration.

%
%
%
%

 
\subsection{}
La première étape de la démonstration consiste à observer que 
les fonctions formelles~$f_i$ sont algébriques sur
le corps~$F(X)$ des fonctions rationnelles sur~$X$.

Pour $i\in I$, posons $Q_i=(P_i,f_i(P_i))$; c'est un point
rationnel de $X\times\P_1$.
Le graphe des fonctions formelles $(f_i)$ fournit
un sous-schéma formel~$\widehat\Gamma$ de $Y=X\times\P_1$ 
supporté par la famille finie $(Q_i)_{i\in I}$.
Par hypothèse, ce sous-schéma formel 
est la fibre générique d'un sous-schéma formel
lisse d'un schéma de la forme $\mathscr X\otimes_R\P_1$,
où $\mathscr X$ est un modèle propre et lisse de~$X$
sur un sous-anneau de type fini de~$F$.
Il est en particulier A-analytique.

Pour toute place~$v$ de~$F$, la fonction méromorphe~$f_v$
définit un morphisme de $\Omega_v$ dans $Y_v$ 
qui induit, au voisinage de tout point~$z_i$,
un isomorphisme du complété formel de~$\Omega_v$ en~$z_i$
sur le germe formel~$\widehat\Gamma_{Q_i,F_v}$.

Notant~$C$ l'adhérence de Zariski de~$\widehat\Gamma$ dans~$Y$,
il résulte donc du théorème~\ref{theo.algebrique} que $C$ est une courbe.

\subsection{}\label{ss.intro-p}
Soit $X'$ la normalisée de~$C$ et 
soit $p\colon X' \to X$ la composition du morphisme
de normalisation et de la première projection.
Puisque le schéma formel $\widehat\Gamma$ est lisse, 
c'est une réunion de branches formelles de $C$.
Il existe ainsi une famille $(P'_i)_{i\in I}$ de points rationnels de~$X'$
telle que le morphisme~$p$ induise, pour tout~$i$,
un isomorphisme du complété formel de~$X'$ en~$P'_i$
sur le complété formel de~$X$ en~$P_i$.
Comme $\Omega_v$ est lisse, la fonction méromorphe~$f_v$ se relève
en un morphisme $f'_v\colon \Omega_v \to X'_v$,
et les fonctions formelles~$f_i$
en des sections formelles $f'_i\colon \widehat X_{P_i}\to \widehat X'_{P'_i}$
de~$p$.

Pour démontrer le théorème~\ref{theo.rationnel},
il suffit de prouver que $p$ est un isomorphisme.
Notons en effet $f\colon X\to\P_1$ 
la composition de l'isomorphisme réciproque~$q$, 
du morphisme de normalisation et de la seconde projection.
Par construction, $f$ se factorise par l'adhérence 
de~$\widehat\Gamma $ dans~$Y$, pour la topologie de Zariski,
de sorte que pour tout~$i$, le complété
formel de~$f$ en~$P_i$ s'identifie à la fonction formelle~$f_i$,
ce qu'il fallait démontrer. 

\subsection{}
Nous avons évoqué avant l'énoncé du théorème~\ref{theo.rationnel}
que la matrice~$\mathmat G$ est associée à un accouplement
en théorie d'Arakelov.  
Donnons maintenant quelques précisions sur sa définition.

Il y a plusieurs façons de construire cet accouplement de hauteurs
qui généralise la théorie d'Arakelov originelle.
Le point de vue que nous adoptons ici requiert la théorie~$W^1$
de~\cite{Bost-1999} aux places archimédiennes;
aux places non archimédiennes, on peut la formuler en termes
de modèles entiers ainsi que nous l'avions fait dans~\citep{BostChambert-Loir-2009} mais je préfère maintenant
la présentation développée par~\citet{Thuillier-2005} 
qui repose sur la théorie de Berkovich.
Nous notons $\hPic(X)$ le groupe abélien 
des fibrés
en droites munis de métriques continues et $W^1$-régulières 
muni de la forme bilinéaire symétrique $\langle\cdot,\cdot\rangle$
tel que définis par~\citet[4.2.2, 4.3]{Thuillier-2005}.

Pour tout~$i$, le fibré en droite $\mathscr O_X(P_i)$
muni des métriques capacitaires données par les ouverts~$\Omega_v$
fournit un élément $\widehat P_i$ de ce groupe abélien
\citep[4.3.11]{Thuillier-2005},
et même de son sous-monoïde défini
par~\citet{Zhang-1995a} des fibrés en droites amples
munis de métriques intégrables.

Pour tout couple~$(i,j)$,
on a alors $\mathmat G_{i,j}=\langle \widehat P_i,\widehat P_j\rangle$.

\begin{prop}\label{prop.markov}
La matrice~$\mathmat G$ est irréductible au sens de la théorie
des chaînes de Markov.
\end{prop}
\begin{proof}
Il s'agit de démontrer que pour tout sous-ensemble non vide~$J$ de~$I$
tel que $I\setminus J$ n'est pas vide, il existe $j\in J$
et $i\in I\setminus J$ tel que $\mathmat G^{i,j}>0$.

Le lemme~\ref{lemm.minimax} fournit une famille $(a_i)_{i\in I}$
d'entiers strictement positifs tels que
$\langle \sum_i a_i \widehat P_i,\widehat P_j\rangle>0$
pour tout $j\in I$.
Posons $\widehat D=\sum_i a_i \widehat P_i$.
C'est un élément \emph{big} et \emph{nef} de $\hPic(X)$.

Posons $\widehat D_1=\sum_{i\in I\setminus J} a_i \widehat P_i$
et $\widehat D_2=\sum_{i\in J} a_i \widehat P_i$.
Puisque \[ \langle \widehat D_1,\widehat D_2\rangle
 = \sum_{i\in I\setminus J} \sum_{j\in J} a_i a_j \mathmat G^{i,j}, \]
il suffit de prouver que $\langle \widehat D_1,\widehat D_2\rangle>0$.
L'assertion à démontrer est ainsi un cas particulier
de l'analogue en théorie d'Arakelov
du lemme~2 de~\cite{Ramanujam-1972}, où l'on remplace
le théorème de l'indice de Hodge pour les surfaces par sa variante
en théorie d'Arakelov.

Le lemme~6.11 de~\citep{Bost-1999} entraîne que
l'on a $\langle \widehat D_1,\widehat D_2\rangle\geq 0$,
et que que si $\langle\widehat D_1,\widehat D_2\rangle=0$,
alors, dans $\hPic(X)$, on a
soit $\widehat D_1=0$, soit $\widehat D_2=0$,
soit il existe $\lambda>0$ tel que $\widehat D_2=\lambda\widehat D_1$.
Les deux premiers cas sont impossibles 
car $\langle\widehat D,\widehat D_1\rangle>0$ et $\langle\widehat D,\widehat D_2\rangle>0$ par la définition de~$\widehat D$.
Dans le dernier cas, on a 
\[ \langle\widehat D,\widehat D\rangle
= \langle (1+\lambda) \widehat D_1, (\lambda^{-1}+1)\widehat D_2\rangle=0, \]
ce qui est également absurde.
Ceci conclut la démonstration.
\end{proof}

\begin{coro}
La courbe~$X'$ est géométriquement intègre.
\end{coro}
\begin{proof}
Puisqu'elle est lisse et qu'elle contient un point rationnel,
il suffit de prouver qu'elle est irréductible. 
D'après la proposition précédente, la matrice~$\mathmat G$
est irréductible.
La remarque~\ref{remas.algebrique}, (4), entraîne donc que
la courbe~$C$ est irréductible. Il en est donc de même de~$X'$.
\end{proof}

\subsection{}
Nous pouvons maintenant conclure la démonstration
du théorème~\ref{theo.rationnel}.
Comme expliqué au~\S\ref{ss.intro-p}, il reste
à démontrer que le morphisme~$p\colon X'\to X$ est un
isomorphisme.

De la stricte positivité
de la valeur~$V_{\mathmat G}$ du jeu de matrice~$\mathmat G$,
Nous avons déduit au cours de la démonstration
de la proposition~\ref{prop.markov}, 
qu'il existe 
une famille $(a_i)_{i\in I}$ d'entiers strictement positifs
tels que $\langle\sum_i a_i \widehat P_i, \widehat P_j\rangle >0$
pour tout~$j\in I$.
L'élément $\widehat D=\sum_i a_i \widehat P_i$ de $\hPic(X)$
est alors numériquement effectif 
et vérifie $\langle \widehat D, \widehat D\rangle>0$.
Puisque la courbe~$X'$ est géométriquement intègre,
la proposition~7.5 de \cite{BostChambert-Loir-2009}
entraîne alors que $p$ est un isomorphisme, comme
il fallait démontrer.

%
%
%

\bibliography{Borel-Dwork}

\end{document}